\documentclass{amsart}
\usepackage{graphicx,euscript,eufrak}

\newtheorem{theorem}{Theorem}[section]
\newtheorem{lemma}[theorem]{Lemma}

\newtheorem{conjecture}[theorem]{Conjecture}

\theoremstyle{definition}
\newtheorem{remark}{Remark}

\def\vc#1{\mathbf #1}
\def\mymatrix#1{\begin{bmatrix}#1\end{bmatrix}}
\def\O{{\mathcal O}}

\long\def\ignore#1{}

\def\H{{\mathcal H}}

\def\ss{\sigma}

\def\Pic{\hbox{Pic}}
\def\Aut{\hbox{Aut}}

\def\A{{\mathcal A}}
\def\Bbb#1{{\mathbb #1}}
\def\PPP{\Bbb P^1\times\Bbb P^1\times \Bbb P^1}

\def\D{{\mathcal D}}
\def\E{{\mathcal E}}

\def\dd{\delta}
\def\V{{\mathcal V}}

\def\M{{\mathcal M}}

\def\K{{\mathcal K}}

\def\qed{\tiny $\clubsuit$ \normalsize}
\def\ord{\mbox{ord}}
\def\C{\mathbb C}
\def\Q{\mathbb Q}
\def\P{\mathbb P}

\begin{document}
\title{K3 surfaces, rational curves, and rational points}
\author{Arthur Baragar}
\author{David McKinnon}
\address{Department of Mathematical Sciences, University of Nevada
 Las Vegas, Box 454020, 4505 Maryland Parkway, Las Vegas, NV 89154-4020, U.S.A.}
  \email{baragar@unlv.nevada.edu}
 \address{Department of Pure Mathematics, University of Waterloo, Waterloo, ON, N2L 3G1, CANADA}

 \email{dmckinnon@uwaterloo.ca}

\begin{abstract}
We prove that for any of a wide class of elliptic surfaces $X$ defined
over a number field $k$, if there is an algebraic point on $X$ that lies on 
only finitely many rational curves, then there is an algebraic point on $X$
that lies on no rational curves.  In particular, our theorem applies to a
large class of elliptic $K3$ surfaces, which relates to a question posed by 
Bogomolov in 1981.  We apply our results to construct an explicit algebraic
point on a $K3$ surface that does not lie on any smooth rational curves.
\end{abstract}
\maketitle


\section{Introduction}

In 1981, Fedor Bogomolov made the following conjecture (\cite{BT}):

\begin{conjecture}\label{bogomolov}
Let $k$ be either a finite field or a number field.  Let $X$ be a K3
surface defined over $k$.  Then every $\overline{k}$-rational point on
$X$ lies on some rational curve $C\subset X$, defined over
$\overline{k}$.
\end{conjecture}

In the number field case, supporting evidence for this conjecture has
been less forthcoming than in the finite field case.  Indeed, in
\cite{BT}, Bogomolov and Tschinkel describe the conjecture as an
``extremal statement'' that is ``still a logical possibility''.  The
purpose of this paper is to provide more evidence that in the number
field case at least, Conjecture~\ref{bogomolov} is false.  (In the
finite field case, progress has been made towards a proof of
Conjecture~\ref{bogomolov} --- in particular, in \cite{BT}, the
authors prove the conjecture for Kummer surfaces defined over a finite
field.)

In particular, the main theorem (Theorem~\ref{main}) proves that for a
wide class of elliptic surfaces $X$, if there is an algebraic point on
$X$ that lies on only finitely many rational curves, then there is an
algebraic point on $X$ that lies on no rational curves.  The same
techniques prove an alternative version of the theorem
(Theorem~\ref{alternative}) that, with an additional hypothesis on the
point, would disprove Conjecture~\ref{bogomolov}.  Unfortunately, it
is not clear to the authors how to construct a specific $K3$ surface
and point $P$ that satisfy the hypotheses of Theorem~\ref{main} or
Theorem~\ref{alternative}.

On the other hand, in section~\ref{surface}, we are able to use these
techniques to find an explicit algebraic point (indeed, four points)
on a K3 surface $X$ that do not lie on any {\em smooth} rational
curve on $X$.  The surface $X$ we use has infinitely many smooth
rational curves, and although it is not difficult to show that there
must exist some algebraic point that does not lie on a smooth rational
curve, the authors are not aware of any explicit constructions of such
a point.  Moreover, our techniques allow, at least in principle, for
any positive integer $d$, the explicit construction of an algebraic
point on a K3 surface $X$ which does not lie on any rational curves of
arithmetic genus at most $d$.  As $d$ increases, however, the calculations
involved quickly become unmanageable, which is why we restrict ourselves 
to the case $d=0$ in the present manuscript.

\section{Main Theorem}\label{maintheorem}

Before we state the main theorem, we will review some elementary
definitions and results about places of curves.

Let $C$ be an irreducible curve defined over the field $\C$ of complex
numbers.  A place of $C$ is a closed point of the normalization
$\tilde{C}$.  A map of curves $h\colon C_1\to C_2$ is ramified at a
place $Q$ to degree $d$ if and only if the corresponding map
$\tilde{h}\colon\tilde{C_1}\to\tilde{C_2}$ is ramified at $Q$ to
degree $d$.

Assume further that $C$ lies on a smooth surface $X$, and let $D$ be a
divisor on $X$ such that $C$ is not contained in the support of $D$.
The intersection multiplicity of $C$ and $D$ at a place $Q$ of $C$ is
$\ord_Q(f^*D)$, where $\ord_Q$ is the discrete valuation associated to
the smooth point $Q$ on $\tilde{C}$, and $f\colon\tilde{C}\to C$ is
the normalization map.  Note that by Example~7.1.17 of \cite{Fu}, we have
the relation:
\[I(f(Q),C\cdot D;X)=\sum_{f(R)=f(Q)}\ord_R(f^*D)\]
In other words, the intersection multiplicity of $C$ and $D$ at a
point $P$ equals the sum of the intersection multiplicities of C with
D at all the places of $C$ lying over $P$.

\begin{theorem}\label{main}
Let $\pi\colon X\to\P^1$ be a smooth elliptic surface defined over a
number field $k$, with a section $S$ and at least five nodal singular
fibres.  Let $\E$ be the elliptic curve over $k(T)$ corresponding to
this fibration with zero section $S$.  Assume that the primitive
2-torsion on $\E$ corresponds to an irreducible curve of positive
(geometric) genus on $X$.

Let $P$ be an algebraic point on $X$, lying on a smooth fibre $E$ of
$\pi$.  Let $[2]\colon E\to E$ be multiplication by 2.  Assume that
there are only finitely many rational curves on $X$ through $P$.  Then
there is a point $Q$ on $X$ such that $[2]^nQ=P$ for some positive
integer $n$ and $Q$ lies on no rational curves on $X$.  In particular,
$Q$ is an algebraic point on $X$.
\end{theorem}

\vspace{.1in}

\noindent
{\it Proof:} \/ Let $f\colon X\to X$ be the rational map corresponding
to the multiplication by 2 on the elliptic curve $\E$.  Then $f$ is
well defined at every point of $X$ except the singular points of the
singular fibres of $\pi$.  Let $m\colon Y\to X$ be a minimal
resolution of $f$ --- that is, assume that $Y$ is a minimal blowup of
$X$ with the property that $f$ extends to a morphism $m\colon Y\to X$.
Let $\psi\colon Y\to X$ be the blowing down map.

It is a straightforward calculation that near a point $Q$ of $X$ that
is the node of a singular fibre of $\pi$, $Y$ is simply the blowup of
$X$ at $Q$.  It is also straightforward to see that $m$ is unramified
over smooth fibres of $\pi$, and that over nodal fibres of $\pi$, $m$
is ramified (to order two) precisely along the curve $\psi^{-1}(Q)$,
where $Q$ is the node.  In particular, $m$ induces an \'etale map of
degree 4 from each smooth fibre of $\pi$ to itself.  Over a nodal
fibre $N$, $m$ induces a map of degree 2 from $N$ to itself, ramified
over the two places lying over the node $Q$, and $m$ restricts to the
normalization map from the curve $\psi^{-1}(Q)$ to $N$, although ---
as previously noted --- $m$ is ramified to order two along
$\psi^{-1}(Q)$.  Thus, in particular, we have the equality of divisors
$m^*N=N+2\psi^{-1}(Q)$.  

The heart of the proof lies in the following lemma:

\begin{lemma}\label{key}
Let $C$ be any irreducible curve on $X$.  If $C$ is not a component of
a singular fibre of $\pi$, then $m^{-1}C$ has at least one component
that is not a rational curve.  Moreover, if $m^{-1}C$ has a component that
is a rational curve, then it is has exactly two components, and $m$ induces
a birational map from the rational component to $C$.
\end{lemma}

\vspace{.1in}

\noindent
{\it Proof:} \/ If $C$ is not itself a rational curve, then clearly
every component of $m^{-1}C$ is not a rational curve.  Thus, we may
assume that $C$ is a rational curve, albeit possibly a singular one.
Since $C$ is not a component of a singular fibre of $\pi$, it follows
that $\pi$ induces a nonconstant morphism $g\colon C\to\P^1$.  Let $d$
be the degree of $\pi|_C$.  That is, let $d=C\cdot F$, where $F$ is
the divisor class of a fibre of $\pi$.  Since $C$ is rational,
Hurwitz's Theorem (\cite{Ha}, Corollary~IV.2.4) implies that $g$ has
ramification degree $2d-2$.

If $d=1$, then $C$ is a section of $\pi$.  If $C=S$, then by
assumption the divisor $m^{-1}C$ has two components: $S$ and the
irreducible 2-torsion, which is assumed to be non-rational.  If $C$ is
not equal to $S$, then it is a translate of $S$, and thus $m^{-1}C$ is
isomorphic to $m^{-1}S$.  The lemma is therefore true for $d=1$, and we
henceforth assume that $d\geq 2$.

We next deal with the case that $m^{-1}C$ is reducible.  Since $m$ has
degree four, the degrees of the components of $m^{-1}C$ over $C$ must
sum to four.  Since the 2-torsion of $\E/k(T)$ is irreducible of
degree 3 over $\P^1$, it immediately follows that there can be no more
than one component of $m^{-1}C$ of degree 1, and no components of
degree 2.  The only remaining reducible case has one component of
degree 1 and one component of degree 3.  The degree 1 component is
clearly rational, so if the degree 3 component were also rational,
then there would be nontrivial 2-torsion of $\E$ defined over a
rational function field, which is impossible since the primitive
2-torsion is non-rational.

We now restrict to the case that $m^{-1}C$ is irreducible.  For any
place $Q$ of $C$, the ramification degree of $g$ at $Q$ is equal to
the intersection multiplicity of $Q$ with the fibre of $\pi$ through
$Q$.  Over the nodal fibres of $\pi$, these intersection
multiplicities sum to at least $5d$, while the ramification degree of
$g$ is $2d-2$.  Thus, since $d\geq 2$, there are at least 8 places $Q$
of $C$ lying on nodal fibres of $\pi$ such that $g$ is unramified at
$Q$.

Let $Q$ be a place of $C$ lying on a nodal fibre of $\pi$, and such
that $g$ is unramified at $Q$.  The intersection multiplicity of
$C$ with the fibre at $Q$ is one, so $Q$ is a smooth point of the
nodal fibre.  This means that $m^{-1}(Q)$ is a set of exactly three
points of $Y$, exactly one of which --- call it $R$ --- lies on the
ramification locus of $m$.  The point $R$ blows down to the node $T$
of the nodal fibre on which $Q$ lies (that is, $\psi(R)=T$), and the
other two points lie on the smooth part of the fibre.

If $R$ corresponds to more than one place of $m^{-1}C$, then $m^{-1}C$
is singular at $T$, and thus has multiplicity at least two at $T$.
Since the fibre is also singular at $T$, this means that the
intersection multiplicity of the fibre with $m^{-1}C$ along $m^{-1}Q$
is greater than 4, which is clearly impossible.  Thus, $R$ corresponds
to a single place of $m^{-1}C$.

But this means that $m|_{m^{-1}C}\colon m^{-1}C\to C$ is ramified at
the place $R$.  Since there are at least 8 such places, it follows
from Hurwitz's Theorem that the curve $m^{-1}C$ is not rational.  This
concludes the proof of the lemma.  \qed

\vspace{.1in}

We now complete the proof of the theorem.  Assume that there are $r$
rational curves on $X$ through $P$, and let $Z$ be one of them.  Since
$P$ lies on a smooth fibre of $\pi$, Lemma~\ref{key} implies that the
set $m^{-1}Z$ has at least one non-rational component $Y$.  Let $Q$ be
a point in $Y\cap m^{-1}(P)$.  The morphism $m$ induces a function
$\M$ from \{rational curves through $Q$\} to \{rational curves through
$P$\}, and by Lemma~\ref{key}, $\M$ is injective, and not surjective
(because $Z$ is not in the image of $\M$).  Thus, there are strictly
fewer rational curves through $Q$ than through $P$.  By iterating this
procedure at most $r$ times, one obtains a point $Q$ such that
$[2]^nQ=P$ and such that no rational curves on $X$ contain $Q$.  \qed

\vspace{.1in}

It seems highly unlikely that every algebraic point on a K3 surface
lies on infinitely many rational curves.  However, it is easy to
construct examples of algebraic points on K3 surfaces that lie on
infinitely many rational curves.  For example, if $P$ is a point that
is fixed by a rational map $f\colon X\to X$, and if $P$ lies on some
rational curve $C$ that is not a pre-periodic curve of $f$ (that is,
the sequence of curves $\{f^n(C)\}$ is not eventually periodic), then
$P$ obviously lies on the infinite set of rational curves
$\{f^n(C)\}$.  This occurs when, for example, the point $P$ is a point
of intersection of the zero section $S$ of an elliptic fibration and a
rational curve $C$ which is non-torsion.  

However, these examples all admit a number field $k$ over which all
the relevant rational curves are defined.  In
Theorem~\ref{alternative}, we describe a possible means of
circumventing this problem.

\vspace{.1in}

\begin{theorem}\label{alternative}
Let $X$ satisfy the conditions of Theorem~\ref{main}, and let $P$ be
an algebraic point on $X$, defined over a field $k$.  Let $f\colon
X\to X$ be the rational map given by multiplication by two.  Assume
that $f^{-1}(P)$ is irreducible over $k$.  Then for any point $Q$
satisfying $f(Q)=P$ (that is, $2Q=P$), there are no rational curves
$C$ on $X$ through $Q$ such that $f(C)$ is defined over $k$.
\end{theorem}

\noindent
{\it Proof:} \/ Let $C/\overline{k}$ be any irreducible curve through
$Q$ (possibly singular, possibly not defined over $k$), and assume
that the curve $f(C)$ is defined over $k$.  Let $G$ be the curve
$\overline{f(C)}$, and let $D=f^{-1}(G)$.  Every point of $f^{-1}(P)$
lies on some component of $D$, and every component of $D$ passes
through some point of $f^{-1}(P)$.  Since the points of $f^{-1}(P)$
are all Galois conjugates, it follows that the components of $D$ are
all Galois conjugates.  By Lemma~\ref{key}, at least one component of
$D$ is non-rational.  It therefore follows that all the components of
$D$ are non-rational.  In particular, $C$ is not a rational curve.
\qed

\vspace{.1in}

Notice that if we assume further that every rational curve through $P$
is defined over $k$, then Theorem~\ref{alternative} implies that there
are no rational curves through $Q$ at all, providing a counterexample
(indeed four counterexamples) to Conjecture~\ref{bogomolov}.  

\section{An explicit example}\label{surface}

In this section, we will use the techniques of the previous sections
to exhibit a specific example of a K3 surface $X$ and a point $P$ on
$X$ such that $P$ does not lie on any smooth rational curves on $X$.

First, we describe a K3 surface $X$ with the following properties:
\begin{enumerate}
 \item  $X$ has an elliptic fibration with a section,
 \item  $X$ contains an infinite number of $(-2)$-curves, all
 defined over $\Bbb Q$,
 \item  at least five of the singular fibers in the elliptic
 fibration are irreducible curves with a nodal singularity,
 \item the divisor of $2$-torsion points is irreducible, and
 \item $X$ contains a rational point $P$ such that the four solutions
 $Q$ to $[2]Q=P$ are all Galois conjugates, where we use the addition
 on the elliptic fiber $E$ that contains $P$, where the zero element
 is the intersection of $E$ with the given section.
\end{enumerate}

Notice that any of the points $Q$ and the surface $X$ will provide the
specific example we seek.  To see this, notice that by
Theorem~\ref{alternative}, there are no rational curves through $Q$
that are defined over $\Bbb Q$.  Since every smooth rational curve on
$X$ is defined over $\Bbb Q$, we conclude that there are no smooth
rational curves on $X$ through $Q$.

The surface we choose comes from the class of K3 surfaces that are defined by
smooth $(2,2,2)$ forms in $\PPP$, have Picard number four, and
include a line parallel to one of the axes.  Such surfaces are studied
in \cite{Ba1}, from which we borrow several basic results.  The
specific surface we consider is the surface $X$ with affine equation
\begin{align*}
F(x,y,z)=&x^2(y^2+2y^2z+yz+z^2+2y+3z)+x(y^2z^2+3y^2z+2y^2+z) \\
&\qquad +(y^2z^2+3y^2z+2y^2+y+z)=0.
\end{align*}
The surface $X$ includes the line $(x,0,0)$, so has Picard number
at least four.  Since $F(x,1/y,z-1)$ is equivalent, modulo $2$,
to the surfaces described by van Luijk in \cite{B-vL}, it has
Picard number at most $4$, so the Picard number is exactly $4$.

\subsection{The Picard group and a fibration with section}
Let $\Bbb P_i^1$ be the $i$th copy of $\Bbb P^1$ in $\PPP$.  Let
$\pi_i$ be projection onto $\Bbb P_i^1$, and $D_i$ the divisor
class of the curve $\pi_i^*(H)$ for some point $H\in \Bbb P^1$.
Let $D_4$ be the divisor class that contains the curve $(x,0,0)$.
Since this curve is a smooth rational curve, its self
intersection is $-2$ (by the adjunction formula).  Since there is
only one $(-2)$-curve in the class $D_4$, we will sometimes abuse
notation and let $D_4$ represent the curve itself too.  The set
$\D=\{D_1,D_2,D_3,D_4\}$ is a basis of $\Pic(X)$ and the
intersection matrix is
\[
J=[D_i\cdot D_j]=\mymatrix{0&2&2&1\\ 2&0&2&0 \\ 2&2&0&0\\
1&0&0&-2}
\]
(see \cite{Ba1}).  The curves $E=\pi_1^*(H)$ are elliptic curves
(again by the adjunction formula) so generate an elliptic
fibration of $X$. The $(-2)$-curve $D_4$ is a section, since
$D_4\cdot D_1=1$.

\subsection{The group of automorphisms and the $(-2)$-curves}

In this subsection, we show that the set of irreducible $(-2)$-curves
on $X$ are all in the $\Aut(X/\Bbb Q)$-orbit of the $(-2)$-curve $D_4$,
so are all rational.  We will first have to describe $\Aut(X/\Bbb Q)$
(or more precisely, a sufficiently large subgroup of $\Aut(X/\Bbb
Q)$).

Let $p_i$ be projection onto $\Bbb P_j^1\times \Bbb P_k^1$, where
$(i,j,k)$ is a permutation of $(1,2,3)$.  Both $p_2$ and $p_3$
are everywhere double covers, and $p_1$ is a double cover at all
points in $\Bbb P_2^1\times \Bbb P_3^1$ except the point
$(0,0)$.  Where we have a double cover, let us define $P'$ by
$p_i^{-1}(p_i(P))=\{P,P'\}$, and set $\ss_i(P)=P'$.  In
\cite{Ba1}, we describe how to extend $\ss_1$ to points on
$D_4$.  These three maps are in $\Aut(X/\Bbb Q)$.

This next automorphism is a little less obvious than those presented
above.  Given a point $P\in X$, let $E$ be the elliptic curve on $X$
that contains $P$ and is in the divisor class $D_1$.  Let $O_E$ be the
point of intersection of $E$ with the section $D_4$.  Define
$\ss_4(P)=-P$, where $-P$ is the additive inverse of $P$ in the group
on $E$ with zero $O_E$.  Then $\ss_4$ is in $\Aut(X/\Bbb Q)$.  In
other words, $\ss_4$ is the automorphism induced by multiplication by
$-1$ on the elliptic fibration corresponding to $D_1$

Let $\A=\langle \ss_1, \ss_2, \ss_3, \ss_4\rangle$ be the group
generated by the automorphisms $\ss_i$.  To understand how $\A$ acts
on $D_4$, we look at its action on the Picard group.  The main result
of \cite{Ba1} is to describe the ample cone $\K$ for surfaces in a
class of K3 surfaces that contains $X$.  In the basis $\D$, the group
of symmetries of $\K$ is $O''=\langle S, T_1, T_2, T_4\rangle$, where
\begin{alignat*}{3}
U&=\mymatrix{1&0&0&0 \\ 0&0&1&0 \\ 0&1&0&0 \\ 0&0&0&1},
&\qquad T_1&=\mymatrix{-1&0&0&0\\ 2&1&0&0 \\ 2&0&1&0 \\ -1&0&0&1}, \\
T_2&=\mymatrix{1&2&0&0 \\ 0&-1&0&0 \\ 0&2&1&1 \\ 0&0&0&-1},
&\qquad\hbox{and}\qquad T_4&=\mymatrix{1&8&8&0 \\ 0&-1&0&0 \\ 0&0&-1&0 \\
0&4&4&1}.
\end{alignat*}
It is clear that every automorphism of $X$ acts as a symmetry of
$\K$, so the pullback map sends $\A$ into $\O''$.  In \cite{Ba1},
we further show that the set of irreducible $(-2)$-divisors is
exactly the $O''$-orbit of $D_4$ (this is used to find the faces
of $\K$).  Hence, to show that all irreducible $(-2)$-curves on $X$
are rational, it is enough to find a subgroup of $\Aut(X/\Bbb Q)$
that maps onto $\O''$.  But that might be asking for too much.
Instead, we note that $UD_4=D_4$, and that $\langle U\rangle$ is
a normal subgroup of $\O''$.  Hence, it is enough to find a
subgroup of $\Aut(X/\Bbb Q)$ (namely $\A$) that maps onto
$\O''/\langle U\rangle$ (using the pullback map, modulo $U$).

In \cite{Ba1}, we show that $\ss_i^*=T_i$ for $i=1,2,3$, where
$T_3=UT_2U$.  Let $[T_4]\in \O''/\langle U\rangle$ be
the equivalence class $\{T_4,UT_4\}$.

\begin{lemma}  The image of $\ss_4$ in $\O''/\langle U\rangle$ is
$[T_4]$.
\end{lemma}

\proof  It is clear that $\ss_4(E)=E$ and $\ss_4(D_4)=D_4$, so
$\ss_4^*D_1=D_1$ and $\ss_4^*D_4=D_4$.  This gives us two
eigenvectors of $\ss_4^*$.

Since the intersection pairing defines a Lorentz product (its
signature is $(1,3)$), there is a natural model of hyperbolic
three space in $\Pic(X)\otimes \Bbb R$.  Let $D$ be an ample
divisor (e.g. $D=D_1+D_2+D_3$) and let
\[
\H=\{\vc x\in \Pic(X)\otimes\Bbb R: \vc x\cdot \vc x=D\cdot
D, \vc x\cdot D>0\}.
\]
Define a distance $|AB|$ between points on $\H$ by
\[
(D\cdot D)\cosh(|AB|)=\vc A\cdot \vc B.
\]
Then $\H$ is a model of $\Bbb H^3$.   Since $\ss_4^*$ preserves the
intersection pairing, it is an isometry on $\H$.  Since $\ss_4^*$
fixes $D_1$ and $D_4$, it fixes every point on the line $l$ in $\H$ with
endpoints $D_1$ and $D_1+D_4$.  Note that $\ss_4^2$ is the identity
on $X$.  Thus, $\ss_4^*$ is either the identity, rotation by $\pi$
about the line $l$, or is reflection through a plane that includes
$l$.  The rotation by $\pi$ about $l$ is $T_4$.  Suppose $\ss_4^*$
is reflection through the hyperplane given by $\vc a\cdot \vc x=0$
intersected with $\H$.  Then $\vc a\cdot D_1=0$ and $\vc a\cdot
D_4=0$, so $a_1=-4a_2-4a_3$ and $a_4=-2a_2-2a_3$.  The reflection
$R_{\vc a}$ through $\vc a\cdot \vc x=0$ is given by
\[
R_{\vc a}(\vc x)=\vc x-\frac{2\vc a\cdot \vc x}{\vc a\cdot \vc
a}\vc a.
\]
If $R_{\vc a}\in O''$, then $R_{\vc a}$ must have integer entries.
Since $\vc a$ is the eigenvector of $R_{\vc a}$ with associated
eigenvalue $-1$ (with multiplicity $1$), it can be taken to
have integer entries.  Furthermore, $R_{\vc a}(D_2)$ must have
integer entries.  The second component of $R_{\vc a}(D_2)$ is
\[
\frac{2(a_3^2-a_2^2)}{2a_2^2+3a_2a_3+2a_3^2}.
\]
Let this integer be $k$.  Then
\[
(2k+2)t^2+3kt+(2k-2)=0,
\]
where $t=a_2/a_3$.  Since $t$ is rational, the discriminant
$-5k^2+16$ must be a perfect square, so $k=0$.  Thus $a_2=a_3$ or
$a_2=-a_3$, the first giving $R_{\vc a}=T_4U$ and the second giving
$R_{\vc a}=U$.

Suppose now that $\ss_4^*=I$ or $U$.  Consider the infinite set of
divisors: $C_m=(T_2T_4T_3T_4)^mD_4$.  Since these are all in the
$O''$-orbit of $D_4$, they each represent irreducible $(-2)$-curves,
which we will also denote with $C_m$.  A simple calculation verifies
that $C_m=c_{m,1}D_1+mD_2+mD_3+c_{m,4}D_4$, so $C_m$ is fixed by both
$U$ and $I$.  Another simple calculation verifies that
$T_4T_3T_4T_2D_1=D_1$, so
\[
C_m\cdot D_1=D_4\cdot(T_4T_3T_4T_2)^mD_1=D_4\cdot D_1=1.
\]
Thus, for each fiber $E$, the curve $C_m$ intersects it at exactly
one point, say $P_m$.  The curves $C_m$ for $m=0,...,4$ intersect in
a finite number of points, so there exists a fiber $E$ on which the
five points $P_m$ are distinct.  Since $\ss_4^*(C_m)=C_m$, and there
is only one curve in this class, we get $\ss_4(C_m)=C_m$, so
$\ss_4(P_m)=P_m$.  But by definition, $\ss_4(P_m)=-P_m$, so we get
$2P_m=O$.  Since an elliptic curve has at most four $2$-torsion
points, we arrive at a contradiction.  Thus $\ss_4^*=T_4$ or $UT_4$.
\endproof

Consequently, the set of irreducible $(-2)$-curves on $X$ is the
$\A$-orbit of $D_4$, so all irreducible $(-2)$-curves on $X$ are
defined over $\Bbb Q$.  We will later prove $\ss_4^*=T_4$, though this
refinement is not necessary for our construction.

\subsection{The singular fibers}

The affine singularities on the singular fibers satisfy the
following system of equations
\begin{align*}
F(x,y,z)&=0 \\
\frac{\partial}{\partial y}F(x,y,z)&= 0 \\
\frac{\partial}{\partial z}F(x,y,z)&=0.
\end{align*}
Maple has no problem solving this system of equations, and finds
that $x$ is a root of a polynomial $g(t)\in \Bbb Q[t]$ of degree $24$,
and that $y$ and $z$ are rational functions in $x$.  We check $g(t)$
modulo several different primes, and discover that modulo $13$, $g(t)$
factors into irreducible polynomials of degree $1$ and $23$, with the
root $x=7 \pmod{13}$.  The singularity on this fiber is at $(9,5)
\pmod{13}$, and
\[
F(7,y+9,z+5)=8y^2z^2+8yz^2+8y^2+2yz+6z^2 \pmod{13}.
\]
Since the quadratic part $8y^2+2yz+6z^2$ is irreducible modulo $13$,
the singularity is nodal, so the fiber over the root of $g(t)$ that
reduces to $7$ modulo $13$ is nodal over $\Bbb C$.

We also discover that, modulo $11$, $g(t)$ has no linear factors,
so $g(t)$ is irreducible over $\Bbb Q$ (of course, the rational
root theorem works too). Thus, the singular fibers over each root
of $g(t)$ are all Galois conjugates of each other, so are all
nodal. Furthermore, (though this is not necessary for our
argument), it is well known that an elliptic fibration on a K3 surface
has at most 24 singular fibres (see for example \cite{IS}), so we
have found all of them.  That is, the fiber at infinity is not
singular, and there are no fibers with singularities at infinity.

\subsection{Addition on the fibers and $2$-torsion points}

A fiber $E$ is a $(2,2)$ form, and is a curve of genus $1$.  We
define a `chord and tangent' addition using the intersections of
$(1,1)$ forms with $E$.  Such intersections include four points,
so our definition of addition is a bit tricky.  A $(1,1)$ form is
uniquely defined by three points.  The curve $E$ intersects the
section at one point, which we choose to be $O$. There exists a
$(1,1)$ form that intersects $E$ at $O$ with multiplicity $3$; it
intersects $E$ again at, say, $O'$.  We define $A*B$ to be the
point $C$ such that the $(1,1)$ form through $A$, $B$, and $O'$
intersects $E$ again at $C$.  Then we define $A+B=(A*B)*O$.  It
is useful to observe that $(A*O)*O=A$, and that if $A*B=C$, then
$A*C=B$.

Suppose that, using our definition of addition, $P+Q+R=O$.  Then
\begin{align*}
(((P*Q)*O)*R)*O&=O \\
((P*Q)*O)*R&=O*O=O \\
((P*Q)*O)*O&=R \\
P*Q&=R.
\end{align*}
Thus $P$, $Q$, $R$, and $O'$ all lie on a $(1,1)$ form.  So do $O$
with multiplicity $3$ and $O'$, so as divisors,
$[P]+[Q]+[R]-3[O]=0$.  This shows that our definition is in fact the
usual addition on an elliptic curve.

To solve $[2]P=O$, we note that
\begin{align*}
O=[2]P&=(P*P)*O \\
O*O&=((P*P)*O)*O \\
O&=P*P.
\end{align*}
Thus, we must solve for $P$ such that the $(1,1)$ form through $O$
and $O'$ has a double root.

Let us consider the fiber with $x=0$:
\[
F(0,y,z)=y^2z^2+3y^2z+2y^2+y+z=0.
\]
Our zero is $O=(0,0)$. Let our $(1,1)$ form be (in affine
coordinates) $z=\frac{ay+b}{cy+d}$ and first assume $ad-bc\neq 0$.  Since this
form goes through $O$, we get $b=0$.  When we plug our $(1,1)$ form
into $F(0,y,z)$, we get a factor of $y$ in the numerator. Forcing
$O$ to be a double root, we get $ad+d^2=0$, and since $d\neq 0$,
$a=-d$.  Forcing $O$ to be a triple root, we get $cd=2d^2$, so
$c=-2d$.  Since we could solve for the $(1,1)$ form under the assumption that $ad-bc\neq 0$, and because the $(1,1)$ form through $O$ with multiplicity $3$ is unique, we do not need to consider the cases that correspond to $ad-bc=0$.  Thus, the $(1,1)$ form that intersects $O$ with
multiplicity $3$ is $\frac{y}{2y-1}$, and its fourth point of
intersection is $O'=(7/15,-7)$. The $(1,1)$ forms through $O$ and
$O'$ have the form
\[
z=-\frac{(7c+15d)y}{cy+d}.
\]
Plugging this into $F(0,y,z)$ and dividing through by $y(y-7/15)$,
we get a fraction with numerator
\[
(15d^2+2c^2+11cd)y^2+(c^2+4d^2+4cd)y+2d^2.
\]
This has a double root if its discriminant is zero, which gives us
\[
t^4-36t^2-116t-104=(t+2)(t^3-2t^2-32t-52)=0.
\]
where $t=c/d$.  The solution $t=-2$ gives us the $(1,1)$ form above
that comes from $[2]O=O$, and the other factor is irreducible over
$\Bbb Q$.  Thus the $2$-torsion points on this fiber are not rational.
Hence, they are Galois conjugates of each other.  Thus, the
$2$-torsion divisor on $X$ must be irreducible.

\subsection{An explicit point}
We now describe a point $P$ on a fiber $E$ such that the four
solutions $Q$ to $2Q=P$ are Galois conjugates of each other.  We pick
the fiber $E$ given by $x=0$ and solve for $Q$ such that
$Q*Q=O'=(7/15,-7)\in E$ (or $(0,7/15,-7)$ as a point on $X$).  Thus,
we are solving for $Q$ such that
$2Q=O'*O=\left(0,\frac{-203}{92},\frac{-2198}{841}\right)$.  A $(1,1)$
form through $O'$ with multiplicity two is of the form:
\[
\gamma(y)=\frac{(847c+6525d)y-(1421c+5243d)}{314(cy+d)}.
\]
for some rational numbers $c$ and $d$.  The numerator of
$F\left(0,y,\gamma(y)\right)$ is the product of $(y-7/15)^2$ and a
quadratic.  We let the discriminant of the second quadratic be zero so
that we will have another double root.  This gives us an irreducible
quartic
\[
p(t)=157t^4+2842t^3+19212t^2+57990t+67147,
\]
each of whose roots gives us a distinct $(1,1)$ form.  Let $\zeta$
be one of the roots of $p(t)$.  Solving for where the resulting
$(1,1)$ form intersects the fiber $E$, we obtain the point
\[
Q=\left(0,\frac{1873}{2714}\zeta^3+\frac{1896629}{213049}\zeta^2 +\frac{16345885}{426098}\zeta+\frac{12302005}{213049},\frac{-1}2\zeta^2 - \frac{1421}{314}\zeta-\frac{1758}{157}\right).
\]
The other three solutions to $2Q=O'*O$ are, of course, the Galois
conjugates of $Q$.  By Theorem~\ref{alternative}, these points $Q$ lie
on no rational curves defined over $\Q$, and therefore on no smooth
rational curves.

\subsection{An aside}
We close this section with a proof that $\ss_4^*=T_4$.  As mentioned earlier, this lemma is not necessary for our construction.

\begin{lemma}  The pullback of $\ss_4$ is $T_4$.
\end{lemma}

\proof  The image of $D_4$ under $\ss_3$ is the $(-2)$-curve $D_2-D_4$.  Let $M$ be its image under $\ss_4$.  For a fixed $x$, let $P$ be the unique point of intersection between $D_2-D_4$ and the elliptic curve $E$ over $x$.   We find $-P$ by considering the $(1,1)$ form through $O$, $O'$, and $P$.  Extended over all values of $x$, this gives us a surface $Y$ in $\PPP$; it is an $(r,1,1)$ form for some non-negative integer $r$.  Let $L$ be the curve of points $O'$ as $x$ varies.

We now look at divisors in the space $\PPP$.  Let $B_i=p_i^*(H)$ for a line $H$ in $\Bbb P_j^1\times\Bbb P_k^1$; and let $B_i'=\pi_i^*(H)$ for a point $H$ in $\Bbb P_i^1$.  Then $B_i'\cdot B_j'=B_k$ where $(i,j,k)$ is a permutation of $(1,2,3)$; and $B_i\cdot B_j'=\dd_{ij}$.  The divisor class that contains $X$ is $2B_1'+2B_2'+2B_3'$; the divisor class that contains $Y$ is $rB_1'+B_2'+B_3'$.  The intersection of $X$ and $Y$ is the union of the four curves $D_4$, $D_2-D_4$, $L$ and $M$.  As divisors, $[X]\cdot [Y]=4B_1+(2r+2)B_2+(2r+2)B_3$.  Thus,
\[
[L]+[M]+[D_2]+[D_2-D_4]=4B_1+(2r+2)B_2+(2r+2)B_3.
\]
Since $D_2$ is the intersection of $X$ with $B_2'$, we get $[D_2]=2B_1+2B_3$.  Hence,
\[
[L]+[M]=2B_1+(2r+2)B_2+2rB_3.
\]
By symmetry, $[L]\cdot B_2'=[L]\cdot B_3'$; let this value be $t$, so $[L]=B_1+tB_2+tB_3$.  Then
\[
[M]=B_1+(2r+2-t)B_2+(2r-t)B_3.
\]
Thus,
\[
[M]\cdot (B_2'-B_3')=2.
\]
But $[M]\cdot(B_2'-B_3')=\ss_4^*(D_2-D_4)\cdot (D_2-D_3)$.  If $\ss_4^*=UT_4$, then this last quantity is $-2$, a contradiction.  Thus, $\ss_4^*=T_4$.
\endproof

\begin{remark}  The curve $L$ of points $O'$ is the curve $C_1$ noted earlier.
\end{remark}

Note that this argument can be generalized much further, at least in
principle.  For example, we never used the fact that $(-2)$-curves are
smooth; we only used the fact that they are all defined over $\Q$.  If
we were given a set of rational curves all defined over some number
field $k$ (say, for example, the set of rational curves of arithmetic
genus at most $d$ on a K3 surface), and a $k$-rational point $P$ such
that the divisor $[2]^{-1}P$ is irreducible over $k$, then we would be
able to deduce that any point $Q$ such that $[2]Q=P$ does not lie on
any curve defined over $k$.

The Hilbert Irreducibility Theorem suggests that such points $P$
should be plentiful, given $k$, but computing the field $k$ for large
$d$ is a more daunting task.  One would have to compute a finite set
$\V$ of rational curves such that any rational curve of arithmetic
genus at most $d$ is conjugate to a curve in $\V$ by some automorphism
of $X$, and then compute the splitting field of $\V$.  Since this calculation
likely grows at least exponentially with $d$ (for example, the Yau-Zaslow
conjecture on the number of rational curves in a given divisor class on a 
K3 surface implies this), it seems that our approach is in practice limited
to relatively small $d$.

\ignore{

@article {MR2299785,
    AUTHOR = {Baragar, Arthur and van Luijk, Ronald},
     TITLE = {{$K3$} surfaces with {P}icard number three and canonical
              vector heights},
   JOURNAL = {Math. Comp.},
  FJOURNAL = {Mathematics of Computation},
    VOLUME = {76},
      YEAR = {2007},
    NUMBER = {259},
     PAGES = {1493--1498 (electronic)},
      ISSN = {0025-5718},
     CODEN = {MCMPAF},
   MRCLASS = {14G40 (11G50 14C22 14J28)},
  MRNUMBER = {MR2299785},
}

}


\begin{thebibliography}{B-vL}

\bibitem[Ba1]{Ba1} { A.~Baragar}, ``The ample cone for a K3 surface,''
to appear.

\bibitem[B-vL]{B-vL}  {A.~Baragar, R.~van Luijk}, ``$K3$ surfaces
with Picard number three and canonical vector heights,'' {\it
Math. Comp.}, {\bf 76}(259), 1493 -- 1498 (2007).  MR2299785

\bibitem[Be]{Be} Beauville, A., {\it Complex algebraic
surfaces}, second edition. London Mathematical Society Student Texts,
34. Cambridge University Press, Cambridge, 1996.

\bibitem[BT]{BT} Bogomolov, F.; Tschinkel, Yu., ``Rational curves and
points on $K3$ surfaces'', Amer. J. Math. 127 (2005), no. 4, 825--835.

\bibitem[Fu]{Fu} Fulton, W., {\it Intersection Theory}, second
edition, Springer, New York, 1998.

\bibitem[Ha]{Ha} Hartshorne, R., {\it Algebraic Geometry}, Springer,
New York, 1977.

\bibitem[IS]{IS} Iskovskikh, V. A.; Shafarevich, I. R., {\it Algebraic
surfaces}. Algebraic geometry, II, 127--262, Encyclopaedia Math. Sci.,
35, Springer, Berlin, 1996.

\bibitem[Se]{Se} Serre, J.-P., {\it Lectures on the Mordell-Weil Theorem},
Vieweg, Wiesbaden, 1997.

\end{thebibliography}
\end{document}